\documentclass[11pt]{article}
\usepackage{amsmath,amssymb,amsfonts}

\newtheorem{theorem}{Theorem}

\newtheorem{corollary}{Corollary}

\newtheorem{lemma}{Lemma}

\def\ff#1#2{
  \if#2\empty{\mathbb F}_{#1}
  \else{\mathbb F}_{#1^{#2}}
  \fi}
  
\def\ffs#1#2{
  \if#2\empty{\mathbb F}_{#1}^\ast
  \else{\mathbb F}_{#1^{#2}}^\ast
  \fi}
  
\def\ffx#1#2{
  \if#2\empty{\mathbb F}_{#1}[X]
  \else{\mathbb F}_{#1^{#2}}[X]
  \fi}
  
\def\ffxi#1#2#3{
  \ifnum#3=0{
    \if#2\empty{\mathbb F}_{#1}[X_1,\ldots,X_e]
    \else{\mathbb F}_{#1^{#2}}[X_1,\ldots,X_e]
    \fi}
  \else{
    \ifnum#3=1\ffx{#1}{#2}
    \else{
      \ifnum#3=2{
        \if#2\empty{\mathbb F}_{#1}[X_1,X_2]
        \else{\mathbb F}_{#1^{#2}}[X_1,X_2]
        \fi}
      \else{
        \if#2\empty{\mathbb F}_{#1}[X_1,\ldots,X_#3]
        \else{\mathbb F}_{#1^{#2}}[X_1,\ldots,X_#3]
        \fi}
      \fi}
    \fi}
  \fi}
                                                                                
\newcommand{\lne}{{\mathcal L}}
\newcommand{\Z}{{\mathbb Z}}

\newenvironment{proof}{\noindent {\bf Proof:}}{$\Box$ \vspace{1ex}}

\begin {document}
\title{\bf
On the splitting case of a semi-biplane construction}
\author{Robert S. Coulter
  \thanks{Department of Mathematical Sciences, 
          University of Delaware, Newark, DE, 19716, 
          U.S.A.
A version of this paper appeared in the Proceedings of the 15th AWOCA conference, 2004, pages 192--198, but is not widely available.
     }
\and Marie Henderson}

\date{}
\maketitle
\begin{abstract}
We consider the case where a particular incidence structure
splits into two substructures. The incidence structure in question was
used previously by the authors to construct semi-biplanes $sbp(k^2,k)$
or $sbp(k^2/2,k)$. A complete description of the two substructures is
obtained. We also show that none of the three semi-biplanes,
$sbp(18,6)$, can be described using this construction.
\end{abstract}

\section{Introduction.}
Let $G$ and $H$ be finite abelian groups written additively and of the same
even order $k$.
We call a function $f:G \rightarrow H$ a {\em semi-planar} function if for
every non-identity $a \in G$ the equation
\begin{equation*}
\Delta_{f,a}(x)=f(x+a)-f(x) = y,
\end{equation*}
with $y \in H$, has either 0 or $2$ solutions $x \in G$. Semi-planar
functions are better known in communication security as almost perfect
non-linear functions, see \cite{nyberg92}, or differentially
$2$-uniform functions, see \cite{nyberg93}.  

A {\em semi-biplane}, or $sbp(v,k)$, is a connected incidence
structure which satisfies the following.
\begin{enumerate}

\item Any two points are incident with 0 or $2$ common lines.
\item Any two lines are incident with 0 or $2$ common points.
\end{enumerate}
Such a design contains $v$ points, and $v$ lines with every point
occurring on $k$ lines, and every line containing $k$ points.
In \cite{coulter99a} the authors developed the following method for
constructing semi-biplanes using semi-planar functions. 

Let $G$ and $H$ be as above and let $f:G \rightarrow H$. Define the
incidence structure $S(G,H;f)$ by: 
\begin{align*}
\text{Points: } &(x,y) \text{ with $x \in G$ and $y \in H$} \\
\text{Lines: } &\lne(a,b) \text{ with $a \in G$ and $b \in H$} \\
\text{Incidence: } &(x,y) \ {\rm I} \ \lne(a,b) \Leftrightarrow y=f(x-a)+b.
\end{align*}
When the context is clear, we shall denote the incidence structure
simply by $S(f)$. The following is Proposition 9 of \cite{coulter99a}. 

\begin{lemma}\label{lem1}
Let $G$ and $H$ be finite abelian groups written additively and of the
same even order $k$. Let $f:G\rightarrow H$ be a semi-planar
function. If $S(G,H;f)$ is connected, then it is a $sbp(k^2,k)$. If
$S(G,H;f)$ is not connected, then $S(G,H;f)$ splits into two
sub-structures; both are $sbp(k^2/2,k)$. 
\end{lemma}

So the designs of \cite{coulter99a} are either connected or consist of
two separate substructures of equal size. Also from \cite{coulter99a}
is the following.

\begin{lemma}\label{lem2}
Let $G$ and $H$ be finite abelian groups written additively and of the
same even order $k$. Let $f:G\rightarrow H$ be a semi-planar
function. If $f$ is a bijection, then $S(G,H;f)$ is connected unless
$k=2$. 
\end{lemma}

As there are bijective semi-planar functions known over the additive
group of any finite field $\ff{q}{}$, with $q=2^e$ and $e\ge 3$, it 
follows that there exist $sbp(2^{2e},2^e)$ for all integers $e\ge
3$. There is only one known class of non-bijective semi-planar 
functions: the monomials $f(X)=X^{2^\alpha+1}$ over $\ff{2}{e}$ are
semi-planar if and only if $(\alpha,e)=1$. Here, too, it can be shown
that $S(f)$ is connected provided $e\ge 3$, see Lemma 11 of
\cite{coulter99a}. When $e=2$, then we must have $f(X)=X^3$ and $S(f)$
splits into two identical copies of the hypercube $H(4)$ ($H(k)$ is
the semi-biplane whose incidence graph is the graph of the
$k$-dimensional hypercube). As $H(k)$ is a $sbp(2^{k-1},k)$, it is
easily seen that $S(G,H;f)$ can only describe a hypercube in this case. 

In this paper we are interested in the case where $S(f)$
splits into two substructures, as at this point the only known
examples which do this are the degenerate case where $k=2$ or the
case $k=4$ with $G=H=\Z_2^+\times \Z_2^+$, the hypercube case. We look
at the general theory for the case where $S(f)$ splits in Section 2.
Our main result gives a complete description of the two substructures
in this case, see Theorem \ref{mainthm}. Proposition 16 of
\cite{wild82} shows that there are exactly three non-isomorphic
$sbp(18,6)$, while there are no $sbp(36,6)$. So if a semi-planar
function $f$ over $\Z_6^+$ exists, then $S(\Z_6^+,\Z_6^+;f)$ must split
into two substructures. In Section 3, we show that no semi-planar
function exists over $\Z_6^+$ and hence none of the three $sbp(18,6)$
can be described by the construction of \cite{coulter99a}. 

\section{General Theory}

For each pair $a\in G$, $b\in H$ define 
\begin{equation*}
S(a,b)=\{t\in G\, :\, f(t-a)=f(t)+b\}.
\end{equation*}
Note that if $f$ is semi-planar, then for each pair $(a,b)\in G\times H$ with
$a\ne 0$, either $|S(a,b)|=2$ or $|S(a,b)|=0$.  

\begin{lemma}\label{lem4}
Let $G$ and $H$ be two finite abelian groups (written additively) of
even order $k$ and $f:G\rightarrow H$ be semi-planar. For each pair
$a\in G$, $b\in H$, with $a\neq 0$, $|S(a,b)|=2$ if and only if
\begin{equation*}
\lne(\alpha a, d+b)\cap\lne((\alpha+1)a,d)\neq\varnothing
\end{equation*}
for all $d\in H$ and $\alpha\in\Z$.
\end{lemma}

\begin{proof}
Let $a\in G$, $b\in H$ and $a\neq 0$. For all $d\in H$ and
$\alpha\in\Z$, the lines $\lne(\alpha a, d+b)$ and
$\lne((\alpha+1)a,d)$ intersect (twice) if and only if
$y=f(x-\alpha a)+d+b=f(x-(\alpha+1)a)+d$. Equivalently,
\begin{equation*}
f(x-(\alpha+1)a)-f(x-\alpha a)=b
\end{equation*}
has two solutions. Substituting for $z=x-\alpha a$, we have
$f(z-a)-f(z)=b$ has two solutions, or in other
words, the lines intersect if and only if $|S(a,b)|=2$.   
\end{proof}

A semi-biplane is called {\em divisible} if the points can be partitioned into
classes so that the following property holds: two points from a class
lie on no common line and two points from different classes lie on
exactly two lines.  

\begin{theorem}\label{thm3}
Suppose $G$ and $H$ are two finite abelian groups (written additively)
of even order $k$ and $f:G\rightarrow H$ is semi-planar. If
$S(G,H;f)$ splits into two substructures, then the resulting
$sbp(k^2/2,k)$ are both divisible. 
\end{theorem}

\begin{proof}
A useful property of $S(f)$ is that it is self-dual, see Theorem 7 of 
\cite{coulter99a}. Hence we need only show the equivalent statement
holds for lines. Let $S_1$ and $S_2$ be the two substructures of
$S(G,H;f)$. Let  
\begin{equation*}
P_a=\{b\in H\, :\, \lne(a,b)\in S_1\}
\end{equation*}
for each $a\in G$. We will show that the set $\{P_a\, :\, a\in G\}$ gives the
required classes. From the proof of Proposition 9 of \cite{coulter99a} there
are exactly $k/2$ elements in each set $P_a$. Also, every point of $S_1$ is in 
$\bigcup_{b\in P_a}\lne(a,b)$ as $\lne(a,b_1)\cap\lne(a,b_2)=\varnothing$ for
all distinct $b_1,b_2\in P_a$ and $S_1$ contains exactly $k^2/2$ points.  

Now choose distinct $a,c\in G$. We claim that
$\lne(a,b)\cap\lne(c,d)\neq\varnothing$ for each $b\in P_a$ and 
$d\in P_c$. If this was not the case, then there is a non-empty list
of lines from $P_a$ which have a common point with $\lne(c,d)$, say
$\lne(a,b_1),\lne(a,b_2),\ldots,\lne(a,b_t),$ where $t<k/2$. From the
definition of incidence, and as $f$ is a semi-planar function, we have
a pair of solutions $(x,y)$ for each member of the above list, given by
\begin{equation*}
y=f(x-a)+b_i=f(x-c)+d.
\end{equation*}
By substituting $z=x-a$ we obtain
\begin{equation*}
\Delta_{f,c-a}(z)=f(z-(c-a))-f(z)=b_i-d.
\end{equation*}
In other words, $\Delta_{f,c-a}(z)=b_i-d$ has $2$ solutions $z\in G$ for
each $1\leq i\leq t$. Overall, this accounts for $2t<k$ of the $k$ values
of $\Delta_{f,c-a}(z)$. The remaining values of $\Delta_{f,c-a}(z)$
must therefore correspond to elements $b\in H$ for which $\lne(a,b)$
and $\lne(c,d)$ intersect and $\lne(a,b)\in S_2$. However this
contradicts the assumption that $S(G,H;f)$ splits into two
substructures. It follows that $|\lne(a,b)\cap\lne(c,d)|=2$ for any
$b\in P_a$ and $d\in P_c$ where $a,c\in G$ are distinct while
$\lne(a,b_1)\cap\lne(a,b_2)=\varnothing$ where $b_1,b_2\in H$ are
distinct. Hence $S_1$ is divisible. A similar argument shows $S_2$ is
also divisible. 
\end{proof}

\noindent Note that from the above proof we know that if $S(f)$
splits, then the lines $\lne(a,b)$ and $\lne(c,d)$ from the same
substructure must intersect when $a\neq c$. This will be used
extensively in what follows. 

For the remainder of this section we suppose $f:G\rightarrow H$ is semi-planar,
$|G|=|H|=k>2$, and $S(f)$ splits into two substructures $S_1$ and $S_2$ with
$\lne(0,0)\in S_1$. Note that, by Lemma \ref{lem2}, $f$ is not a bijection.
For $i=1,2$, define 
\begin{equation*}
P_a^i=\{b\in H\, :\, \lne(a,b)\in S_i\}.
\end{equation*}
For each $a\in G$, $P_a^1\cap P_a^2=\varnothing$ while
$P_a^1\cup P_a^2=H$, so the subsets $P_a^1$ and $P_a^2$ of $H$,
partition $H$. 

\begin{lemma}
For non-zero $a\in G$, 
\begin{align*}
P_a^1&=\{b\in H\, :\, |S(a,b)|=2\},\\
P_a^2&=\{b\in H\, :\, |S(a,b)|=0\}.
\end{align*}
\end{lemma}

\begin{proof}
As $P_a^1$ and $P_a^2$ partition $H$ then we need only consider one of
the subsets, say $P_a^1$. By Lemma \ref{lem4},
$\lne(0,0)\cap\lne(a,b)\neq\varnothing$ if $|S(a,b)|=2$. But Theorem
\ref{thm3} shows, by duality, that for $a\neq 0$, $\lne(a,b)\in S_1$ if
and only if $\lne(a,b)$ and $\lne(0,0)$ intersect.
\end{proof}

\begin{theorem}\label{thm6}
The set $P_0^1$ is the subgroup of $H$ of index $2$ and $P_0^2$ is its
coset.  
\end{theorem}

\begin{proof}
If $0\in P_a^1$, then $\lne(0,d)\cap\lne(a,d)\neq\varnothing$ for
all $d\in H$, by Lemma \ref{lem4}. Thus $P_a^i=P_0^i$ in this case. 
Let 
\begin{equation*}
T=\{a\in G\, :\, a\neq 0\, \land\, |S(a,0)|=2\}.
\end{equation*}
As $f$ is not a bijection, there exists a non-zero $a\in G$ for which
$f(t-a)=f(t)$ has a solution, which implies $T$ is non-empty.
Let $a\in T$. For $b_1\in P_a^1$, $|S(a,b_1)|=2$ and $f(x-a)=f(x)+b_1$
has two solutions. Hence for any $b_2\in P_a^1$, we must have
$f(x-a)+b_2=f(x)+b_1+b_2$ has two solutions, or equivalently,
$\lne(a,b_2)\cap\lne(0,b_1+b_2)\neq\varnothing$. As $b_2\in
P_a^1=P_0^1$ and $\lne(a,b_2)$ intersects $\lne(0,b_1+b_2)$, it
follows that $b_1+b_2\in P_a^1=P_0^1$. To summarise, $b_1,b_2\in
P_0^1$ implies that $b_1+b_2\in P_0^1$, that is $P_0^1$ is closed
under addition. It follows that $P_0^1$ is a subgroup of $H$ of
index two. As $P_0^1\cap P_0^2=\varnothing$ while $P_0^1\cup P_0^2=H$,
$P_0^2$ is the coset of $P_0^1$ in $H$. 
\end{proof}

\begin{lemma}\label{lem8}
If $P_a^i\cap P_c^i\neq\varnothing$, for $i=1$ or $i=2$, then
$P_{a-c}^1=P_{c-a}^1=P_0^1$. 
\end{lemma}

\begin{proof}
Suppose $P_a^i\cap P_c^i\neq\varnothing$ where $i=1$ or $i=2$. Then
there exists a $b\in H$ such that
$\lne(a,b)\cap\lne(c,b)\neq\varnothing$. This, in turn, implies that
there is a $t\in G$ for which we have $f(t-a)-f(t-c)=0$. By
substituting for $z=t-c$ we obtain $f(z-(a-c))-f(z)=0$. So
$f(z)=f(z-(a-c))$ which implies
$\lne(0,0)\cap\lne(a-c,0)\neq\varnothing$. As shown in the proof of
Theorem \ref{thm6}, since $0\in P_{a-c}^1$, it follows that
$P_{a-c}^1=P_0^1$ as required. A similar argument shows
$P_{c-a}^1=P_0^1$.  
\end{proof}

Consider the set $A=\{a\in G\, :\, P_0^1=P_a^1\}$.
By Lemma \ref{lem8}, whenever $P_a^i\cap P_c^i\ne\varnothing$ for $i=1$ or
$i=2$, then $a-c\in A$ and $c-a\in A$. Clearly $0\in A$ and $|A|> 1$. For any
$a,c\in A$, successive applications of Lemma \ref{lem8} show $-c\in A$ and
$a-(-c)=a+c\in A$. Hence $A$ is closed and since $G$ is finite, $A$ is a
subgroup of $G$. If $|A|<k/2$, then $|G\setminus A|> k/2$. Now for some fixed
$a\in G\setminus A$ we have 
\begin{equation*}
|\{a-c\, :\, c\in G\setminus A\}|> k/2.
\end{equation*}
But $\{a-c\, :\, c\in G\setminus A\}\subset A$, contradicting
$|A|<k/2$. So we must have $|A|\geq k/2$ and since $A$ is a subgroup
of $G$, $|A|=k/2$ or $A=G$. This proves the following statement,
common in theme with Theorem \ref{thm6}.

\begin{theorem}\label{thm9}
The set $A=\{a\in G\, :\, P_0^1=P_a^1\}$ is either the subgroup of $G$
of index $2$ or $A=G$.
\end{theorem} 

A combination of Theorems \ref{thm6} and \ref{thm9} proves our main
theorem (which shows that if the structure splits, there are only two 
possibilities).

\begin{theorem}\label{mainthm}
Let $f:G\rightarrow H$ be a semi-planar function where $G$ and $H$ are
abelian groups of even order $k$ and $A$ and $B$ the index two subgroups of
$G$ and $H$, respectively. Let $g\in G\setminus A$ and 
$h\in H\setminus B$. If $S(f)$ splits into two substructures $S_1$ and
$S_2$, with $\lne(0,0)\in S_1$, then either 
\begin{enumerate}

\item $\lne(a,b)\in S_1$ if and only if $(a\in G\land b\in B)$, or
\item $\lne(a,b)\in S_1$ if and only if 
$(a\in A\land b\in B)\lor (a\in A+g\land b\in B+h)$.
\end{enumerate}
\end{theorem}

We note that the theorem also holds for the case $k=2$. In this case, $f(x)=x$
or $f(x)=x+1$. In either case, $f$ is a bijection and the splitting structures
correspond to case (ii).
The theorem allows us to show that the two substructures obtained are
isomorphic. 

\begin{corollary}
For any $h\in H\setminus B$, the mapping $\phi_h:G\times H\rightarrow
G\times H$ defined by 
\begin{equation*}
\phi_h(x,y)=(x,y+h)
\end{equation*}
acts as an isomorphism between the two substructures of $S(f)$.
\end{corollary}

Our final general result, which is a simple extension of
\cite{nyberg92}, Proposition 1, will be needed in the next section.

\begin{lemma}\label{lem11}
If $f:G\rightarrow H$ is a semi-planar function, then
\begin{equation*}
\psi(f(\phi(x)+c))+d
\end{equation*}
is a semi-planar function from $G$ to $H$ where $\phi\in Aut(G)$,
$\psi\in Aut(H)$, $c\in G$, and $d\in H$.  
\end{lemma}

\section{The Case $G=H=\Z_6^+$}

In this section we consider the case where $G=H=\Z_k^+$ with $k$ even.
In this case, we represent the mapping $f:\Z_k^+ \rightarrow \Z_k^+$ by
$f=\langle b_0,b_1,\ldots,b_{k-1}\rangle$ where $f(i)=b_i$ for
$0\leq i\leq k-1$.

\begin{lemma}\label{limitlem}
Let $f:\Z_k^+\rightarrow\Z_k^+$ with $k>4$.
If $f(x)=y$ has more than $k/2$ solutions $x\in \Z_k^+$ for any given
$y\in \Z_k^+$, then $f$ is not semi-planar. 
\end{lemma}

\begin{proof}
Suppose that the claim does not hold. Then $f$ is semi-planar and there exists
$y\in\Z_k^+$ such that $|S|>k/2$ where $S=\{x\in \Z_k^+\, :\, f(x)=y\}$.
We wish to show that there exists an $a\in\Z_k^+$ such that $f(x+a)-f(x)=0$ has
more than two solutions. Consider $f=\langle b_0,b_1,\ldots,b_{k-1}\rangle$.
As $|S|>k/2$ there must be two consecutive elements of this list which
are equal. Using Lemma \ref{lem11}, we may assume $b_0=b_1=y$.

If $f$ is semi-planar, then $\Delta_{f,1}(x)=0$ must have two solutions. There
are 2 cases. If $b_2=y$, then we have three consecutive values of $f$ equal to
$y$ and there can be no other consecutive values of $f$ equal. Thus $b_3\ne y$,
and the remaining $k/2-2$ values of $y$ must be placed in $k-4$ places with no
consecutive places equal. It can be seen that the only way to assign the
remaining $y$ values is $b_j=y$ when $j$ is even. Thus, if $k>4$,
$\Delta_{f,2}(x)=0$ has more than two solutions, a contradiction. If
$b_2\ne y$, then there are $k-3$ remaining assignments of which $k/2-1$ must be
$y$ and where $b_{k-1}\ne y$ as this is equivalent to the previous case by
Lemma \ref{lem1}. Provided $k>4$, it follows that $\Delta_{f,2}(x)=0$ has at
least three solutions, contradicting that $f$ is semi-planar.
\end{proof}

It was shown in \cite{wild82} that no $sbp(36,6)$ exists while there
are three non-isomorphic $sbp(18,6)$. It follows that if a semi-planar
function exists over $\Z_6^+$, then the corresponding structure
necessarily splits. We now show that this case is not
possible. Although this might be tested for computationally, a
mathematical proof is preferable.

\begin{theorem}
There is no semi-planar function over $\Z_6^+$.
\end{theorem}

\begin{proof}
Suppose $f$ is a semi-planar function over $\Z_6^+$. By Lemma
\ref{lem11} we may assume that $f(0)=0$ and that no image of $f$
occurs more often than $0\in \Z_6^+$. Further, by Lemma \ref{limitlem},
$f(x)=0$ has at most three solutions. Let 
\begin{equation*}
f=\langle 0,b_1,b_2,b_3,b_4,b_5\rangle.
\end{equation*}
As noted, $S(f)$ must split. As before we denote the two substructures
by $S_1$ and $S_2$ where $\lne(0,0)\in S_1$. It follows from Theorem
\ref{mainthm} that there are two cases.  

First assume $\lne(a,b)\in S_1$ if and only if $b\in \{0,2,4\}$. It
follows that $b_i\in\{0,2,4\}$ and that $|S(a,0)|=2$ for all
$a\in\Z_6^+$. In particular, from $a=1$ there exists two distinct integers
$r,s\in\Z_6^+$ such that $b_{r-1}=b_r$ and $b_{s-1}=b_s$. Appealing to
Lemma \ref{lem11} we may assume, without loss of generality, that
$b_{r-1}=b_r=0$ and $r=1$. Either $b_s=0$ or $b_s\in\{2,4\}$. If
$b_s=0$, then since $f(x)=0$ can have at most three solutions, we must
have $s=2$ and hence $f=\langle 0,0,0,b_3,b_4,b_5\rangle$ with
$b_3,b_4,b_5\in\{2,4\}$. Now $\Delta_{f,3}(\Z_6^+)=\{2,4\}$. However,
as $f$ is semi-planar, the value set of $\Delta_{f,3}$ must have size
three. So $b_s\ne 0$ and $s>2$. As $\phi(x)=-x$ is an automorphism of
$\Z_6^+$, we may assume $b_s=2$ by Lemma \ref{lem11}. There are three
possibilities:  
\begin{align*}
f&=\langle0,0,2,2,b_4,b_5\rangle,\\
f&=\langle0,0,b_2,2,2,b_5\rangle,\\
f&=\langle0,0,b_2,b_3,2,2\rangle.
\end{align*}
In the first case, $b_5\neq 0$, $b_4\neq 2$ and $b_4\neq b_5$. Hence
$b_4=0$. But then $b_5\in\{2,4\}$ and either leads to
$\Delta_{f,2}(x)=2$ having three solutions. Similar arguments remove
the other two possibilities. It follows that no semi-planar function
exists in this case.

Now assume that $\lne(a,b)\in S_1$ if and only if $a,b\in \{0,2,4\}$
or $a,b\in\{1,3,5\}$. This time we have $b_i\equiv i\bmod 2$. By
considering $\Delta_{f,2}(x)$, an application of Lemma \ref{lem11}
shows we may assume that $b_0=b_2=0$ and $b_4=2$. Likewise, we must
have $b_i=b_j$ for a pair $i,j\in\{1,3,5\}$. We first consider the
situation $f=\langle 0,t,0,t,2,v\rangle$ with $t\neq v$. It is
immediate that $t=5$ as otherwise $\Delta_{f,1}(x)=t$ has at least
three solutions. But if $t=5$ then obviously $v\neq 5$, and also, by
considering $\Delta_{f,1}(x)$, $v\neq 1$. So now $t=5$ and
$v=3$. But then $\Delta_{f,3}(x)=3$ has four solutions. It remains to
deal with the case $f=\langle 0,t,0,v,2,v\rangle$. By considering
$\Delta_{f,3}$, it follows that $t=5$ and $v=1$. But then
$\Delta_{f,1}(x)=1$ has three solutions. Hence no semi-planar function
exists in this case either. All possibilities have been exhausted and
the result follows.
\end{proof}

Our last result shows that the splitting case cannot occur when $k=6$.
It is an open problem to determine a semi-planar function over any abelian
group of order $k>4$ where the splitting case occurs.
We conjecture that no such function exists.


\end{document}